\documentclass[12pt, a4paper, oneside]{amsart}
\usepackage{amssymb, amsmath, amsthm, verbatim, cite, graphicx,array}
\usepackage[cp1251]{inputenc}
\usepackage[english]{babel}
\usepackage{longtable, lipsum}

\newcommand{\mypmod}{\hspace{-1mm}\pmod}

\newtheorem{theorem}{Theorem}

\newtheorem{lemma}{Lemma}[section]
\newtheorem{defin}[theorem]{Definition}
\newtheorem{prop}[theorem]{Proposition}
\theoremstyle{definition}
\newtheorem*{prf}{Proof}

\theoremstyle{definition}

\newcommand{\Aut}{\operatorname{Aut}}
\newcommand{\Out}{\operatorname{Out}}
\newcommand{\Inndiag}{\operatorname{Inndiag}}
\newcommand{\diag}{\operatorname{diag}}
\numberwithin{equation}{section}

\begin{document}
\title[The spectra of almost simple groups]{The spectra of almost simple groups with socle $E_7(q)$}

\author{\sc Alexander A. Buturlakin}
\address{Sobolev Institute of Mathematics\\ Koptyuga 4, Novosibirsk 630090, Russia}
\email{buturlakin@gmail.com}

\author{\sc Maria A. Grechkoseeva}
\address{Sobolev Institute of Mathematics\\ Koptyuga 4, Novosibirsk 630090, Russia}
\email{grechkoseeva@gmail.com}

\thanks{Supported by RAS Fundamental Research Program, project  FWNF-2022-0002}

\begin{abstract}
We give an explicit description of the set of element orders for every almost simple group with socle $E_7(q)$.

\end{abstract}

\begingroup
\let\MakeUppercase\relax 
\maketitle
\endgroup

\section{Introduction}

Given a finite group $G$ and a subset $A$ of $G$, we denote the set of the orders of elements of $A$ by $\omega(A)$, and refer to the set $\omega(G)$ as the spectrum of $G$. In this paper we continue to describe the spectra of finite almost simple groups of Lie type, that is, groups $G$ satisfying $L\leq G\leq \Aut L$ for some finite simple group $L$ of Lie type.

The spectra of simple groups of Lie type themselves are known: see \cite{08But.t} for the spectra of linear and unitary groups, \cite{10But.t} for the spectra of symplectic and orthogonal groups, and \cite{18But.t} for the references concerning exceptional groups.

By Steinberg's theorem, the automorphism group $\Aut L$ of a simple group $L$ of Lie type is a split extension of the group $\Inndiag L$ of inner-diagonal automorphisms by a group generated by field and graph automorphisms. A corollary of the Lang--Steinberg theorem, sometimes referred to as the Shintani correspondence, 
allows us to write the set $\omega(L\rtimes \langle \varphi\rangle)$, where $\varphi$ is a field or a graph-field automorphism of $L$, in terms of the spectra of some simple groups (see, for example, \cite[Corollary 14]{06Zav.t}). So the main difficulty is to handle diagonal and graph automorphisms.

The case when $L<G\leq \Inndiag L$ and $L$ is one of $L_n(q)$, $U_n(q)$,  $E_6(q)$, or $^2E_6(q)$ is studied in \cite{22ButGre}. If $L=O_{2n+1}(q)$ with $q$ odd, then $\Inndiag L\simeq SO_{2n+1}(q)$, and the spectrum of $SO_{2n+1}(q)$ can be found in \cite[Corollary 5]{10But.t}.   Also $\omega(G)$ is known when $G$ is the extension of $L$ by a graph automorphism of order $2$ and $L$ is one of $L_n(q)$, $U_n(q)$, $O_{2n}^\pm(q)$ with $q$ odd (see \cite[Lemmas 4.6 and 4.7]{17Gre} for linear and unitary groups and \cite[Lemma 3.7]{18Gr.t} for orthogonal groups).



In the present paper we are concerned with almost simple groups with socle $E_7(q)$. Let $L=E_7(q)$, where $q$ is a power of a prime $p$, and let $x\mapsto \widehat x$ be the natural homomorphism from $\Aut L$ to $\Out L=\Aut L/L$. It is well known that \begin{equation}\label{e:out} \Out L=\langle \widehat \delta\rangle\times \langle \widehat \varphi\rangle,\end{equation} where $\delta=1$ if $q$ is even, $\delta\in \Inndiag L\setminus L$ and $|\widehat \delta|=2$ if $q$ is odd, and $\varphi$ is the field automorphism of $L$  mapping a root element $x_r(s)$ to $x_r(s^p)$.

Our first result, Theorem~\ref{t:E7}, describes the set
$\omega(\Inndiag L)$ for odd $q$. Moreover, it determines the elements of $\omega(\delta L)=\omega(\Inndiag L\setminus L)$ that are maximal under divisibility.
In Proposition~\ref{p:f}, we express $\omega(\alpha L)$ for arbitrary $\alpha\in \Aut L$ in terms of $\omega(L_0)$ and $\omega(\Inndiag L_0\setminus L_0)$, where $L_0=E_7(q_0)$ for some $q_0$. Finally, Theorem \ref{t:main}, a direct consequence of Theorem \ref{t:E7}, Proposition \ref{p:f} and \cite{16But.t}, gives a description of $\omega(G)$ for all almost simple groups $G$ with socle $L$.

Before stating results, we need to introduce some notation. Let $k$ and $m$ be positive integers and let $A$ be a set of positive integers. The greatest common divisor of $k$ and $m$ is denoted by $(k,m)$, while the least common multiple is denoted by $[k,m]$. Given a prime $p$,  we  define $n_p(k)$ to be the minimal power $p^t$ of $p$ such that $p^t>k$.
By $k\cdot A$ we mean the set $\{ka\mid a\in A\}$. We write $\mu(A)$ for the subset of $A$ consisting of numbers maximal under divisibility.
The symbols $\pm$ and $\mp$ are treated as follows. If a formula includes both  $\pm$ and $\mp$, then the corresponding signs should be chosen to be opposite, so, for example, $\{(q^3\pm 1)(q\mp 1)\}$ is a short for $\{(q^3+1)(q-1), (q^3-1)(q+1)\}$. But if it includes $\pm$ twice, the  corresponding signs should be chosen to be equal.

\begin{theorem}\label{t:E7} Let $L=E_7(q)$, where $q$ is a power of an odd prime $p$,  and $G=\Inndiag L$. Let $\nu(q)$ be the union of the following sets:
\begin{enumerate}

\item $\left\{(q^6\pm q^3+1)(q\mp 1)\vphantom{\frac{q^8-1}{2(q\pm 1)}}\right.$, $q^7\pm 1$, $(q^4- q^2+1)(q^3\pm1)$, $(q^5\pm 1)(q^2\mp q+1)$, $(q^5\pm 1)(q\mp1)$,
$\frac{q^8-1}{2(q\pm 1)}$, $\left. q^6-1\vphantom{\frac{q^8-1}{2(q\pm 1)}}\right\}$;

\item $p \cdot \left\{q^5\pm 1,   (q^4+1)(q^2\pm 1),\frac{q^6-1}{2}, (q^3\pm 1)(q^2+1)(q\mp 1), q^4-q^2+1\right \}$;

\item $n_p(3) \cdot \left\{(q^3\pm 1)(q\mp 1), \frac{q^4-1}{2}, \right\}$;

\item $n_p(5)\cdot \{q^3\pm1, (q^2+1)(q\pm1)\}$;

\item $n_p(7) \cdot\{q^2-1\}$;

\item $n_p(11)\cdot \{q\pm 1\}$;
\item $\{n_p(17)\}$,
\end{enumerate}
and let $\nu_{\delta}(q)=\nu(q)\setminus \{p(q^4-q^2+1), n_p(17)\}$.
Then $\omega(G)$ is precisely the set of all divisors of elements of $\nu(q)$.
Furthermore, $\mu(\omega(G\setminus L))\subseteq \nu_{\delta}(q)\subseteq \omega(G\setminus L)$.
\end{theorem}

Recall from \eqref{e:out} that $\varphi$ is a field automorphism of $L$ of order $m$, where $m$ is defined by $q=p^m$. So if $\psi\in\langle \varphi\rangle$, then $|\psi|$ divides $m$ and $q^{1/|\psi|}$ is a positive integral power of $p$.

\begin{prop}\label{p:f} Let $L=E_7(q)$, $\psi\in\langle \varphi\rangle$, $q_0=q^{1/|\psi|}$, and $L_0=E_7(q_0)$. Then $\omega(\psi L)=|\psi|\cdot \omega(L_0)$. If, in addition, $q$ is odd, then $\omega(\psi (\Inndiag L\setminus L))=|\psi|\cdot \omega(\Inndiag L_0\setminus L_0)$.
\end{prop}

In the following definition, the set $\nu_\alpha(q)$ is in fact depends only on $\widehat \alpha$ but it is more convenient to write  $\nu_\alpha(q)$ instead of $\nu_{\widehat \alpha}(q)$.

\begin{defin}\label{d:main} Let $L=E_7(q)$. For every $\alpha\in\Aut L$, define the set $\nu_{\alpha}(q)$ as follows:
\begin{enumerate}
 \item $\nu_1(q)$ is the union of the following sets:

 \begin{enumerate}
 \item $\left\{\frac{(q^6\pm q^3+1)(q\mp1)}{(2,q-1)}\right.$, $\frac{q^7\pm1}{(2,q-1)}$, $\frac{(q^4-q^2+1)(q^3\pm1)}{(2,q-1)}$, $\frac{(q^5\pm1)(q^2\mp q+1)}{(2,q-1)}$, $(q^5\pm1)(q\mp1)$, $\frac{q^8-1}{(q\pm1)(4,q\pm1)}$, $(q^4+1)(q^2-1)$, $q^6-1$, $\left.(q^3\pm1)(q^2+1)(q\mp1)\vphantom{\frac{(q^3\pm1)}{(2,q-1)}}\right\}$;

\item $p \cdot \left\{\vphantom{\frac{(q^3\pm1)}{(2,q-1)}}q^4-q^2+1\right.$, $q^5\pm1$, $\frac{(q^4+1)(q^2\pm 1)}{(2,q-1)}$, $\frac{q^6-1}{(2,q-1)}$, $\left.\frac{(q^3\pm1)(q^2+1)(q\mp1)}{(2,q-1)}\right\}$;

\item $n_p(2) \cdot\left\{\frac{q^5\pm1}{(2,q-1)}\right.$, $\frac{q^6-1}{(q\pm1)(2,q-1)}$, $(q^3\pm1)(q\mp1)$, $\left.q^4-1\vphantom{\frac{(q^3\pm1)}{(2,q-1)}}\right\}$;

\item $n_p(3) \cdot\left\{\frac{(q^3\pm1)(q\mp1)}{(2,q-1)}\right.$, $\left.\frac{q^4-1}{(2,q-1)}\right\}$;

\item $n_p(5) \cdot\left\{\frac{q^3\pm1}{(2,q-1)}, \frac{(q^2+1)(q\pm1)}{(2,q-1)},  q^2-1\right\}$;

\item $n_p(7) \cdot\left\{\frac{q^2-1}{(2,q-1)}\right\}$;

\item $n_p(9)\cdot\left\{q\pm1\right\}$;

\item $n_p(11)\cdot\left\{\frac{q\pm 1}{(2,q-1)}\right\}$;

\item $\{n_p(17)\}$;
\end{enumerate}

 \item if $\psi \in \langle \varphi\rangle$, then $\nu_{\psi}(q)=|\psi|\cdot \nu_1(q^{1/|\psi|})$;
 \item if $q$ is odd, then $\nu_{\delta}(q)$ is as in Theorem \ref{t:E7};
 \item if $q$ is odd and $\psi \in \langle \varphi\rangle$, then $\nu_{\psi\delta}(q)=|\psi|\cdot \nu_{\delta}(q^{1/|\psi|})$;
 \item if $\widehat\alpha=\widehat\beta$, then $\nu_\alpha(q)=\nu_\beta(q)$.
\end{enumerate}
\end{defin}

\begin{theorem}\label{t:main} Let $L=E_7(q)$, $\alpha\in \Aut L$ and let $\nu_\alpha(q)$ be as in Definition \ref{d:main}. Then $$\mu(\omega(\alpha L))\subseteq \nu_{\alpha}(q)\subseteq \omega(\alpha L).$$ In particular, if $L\leq G\leq \Aut L$, then $\omega(G)$ is precisely the set of all divisors of elements of $\cup_{\widehat \alpha\in G/L}\nu_\alpha(q).$
\end{theorem}

Our proof of Theorem \ref{t:E7} is based on Carter's description \cite{78Car} of connected centralizers of semisimple elements and similar to that of the results in  \cite{16But.t}. Since the group $\Inndiag L$  and the universal version of $L$ have isomorphic maximal tori by \cite[Chapter 4]{85Car}, the orders of semisimple elements in Theorem~\ref{t:E7} can be found from  the results of \cite{91DerFak}. The other orders require calculations, some of which are carried out in \textsc{Magma} \cite{MagmaBook}. The proof of Proposition \ref{p:f} is a slight generalization of that of \cite[Proposition 13]{06Zav.t}.

\section{Preliminaries}\label{s:pre}

Let $r$ be a prime. If $a$ is an integer, then $(a)_r$ is the $r$-part of $a$, that, is the highest power of $r$ dividing $a$, and $(a)_{r'}$ is the ratio $|a|/(a)_r$.

Given a finite group $G$, we write $\exp(G)$ for the exponent of $G$, that is, the smallest positive integer $k$ such that $g^k=1$ for all $g\in G$, and $\exp_r(G)$ for the exponent of a Sylow $r$-subgroup of $G$.

As we mentioned in the introduction, our proof is based on Carter's work \cite{78Car}, and we proceed with stating necessary facts and results of this work.

Let $\overline F$ be the algebraic closure of the field of prime order $p$. Let   $\overline G$ be a connected reductive algebraic group over $\overline F$
and $\sigma$ a Steinberg endomorphism of $\overline G$ (that is, $\sigma$ is a surjective endomorphism of $\overline G$ and $C_{\overline{G}}(\sigma)$ is finite). The finite group
$C_{\overline{G}}(\sigma)$ is denoted by $\overline G_{\sigma}$. Also $\overline T$ is a $\sigma$-stable maximal torus of $\overline G$, $W=N_{\overline G}(\overline T)/\overline T$ is the Weyl group of $\overline G$ and $\pi$ is the natural
homomorphism from $N_{\overline G}(\overline T)$ to $W$.

Let $g\in\overline G_\sigma$ and let $s$ and $u$ be the semisimple and unipotent parts of $g$ respectively. Then $g$ lies in the connected component $C_{\overline G}(s)^0$ and $C_{\overline G}(s)^0$ is a $\sigma$-stable reductive subgroup of $\overline G$ of maximal rank. Let $\overline R=C_{\overline G}(s)^0$. We refer to $R=\overline R_\sigma$ as a reductive subgroup of $\overline G_{\sigma}$ of maximal rank. It is clear that $|u|$ divides the $p$-exponent of $R$ and $|s|$ divides the exponent of the center of this group. Thus $|g|$ divides the number $$\eta(R)=\exp_p(R)\cdot \exp(Z(R)).$$ Conversely, it is easy to see that $\eta(R)$ lies in $\omega(\overline G_\sigma)$. So to find $\omega(\overline G_\sigma)$, we need to find the set $\{\eta(R)\}$ where $R$ runs over all reductive subgroups of maximal rank of $\overline G_\sigma$.

By \cite[Section 1.2]{78Car}, every $\sigma$-stable reductive subgroup of $\overline G$ of maximal rank is equal to $\overline R^g$, where $\overline R$ is a $\sigma$-stable reductive subgroup of $\overline G$ containing $\overline T$, $g\in N_{\overline G}(\overline T)$ and $g^\sigma g^{-1}\in N_{\overline G}(\overline T)\cap N_{\overline G}(\overline R)$. In turn, a $\sigma$-stable reductive subgroup of $\overline G$ containing $\overline T$ is determined by its root system, which is a closed subsystem of $\Phi$.

In the remaining part of this section, let $\overline R$ be a $\sigma$-stable reductive subgroup of $\overline G$ containing $\overline T$, $\Phi_1$ the root system of $\overline R$ (with respect to $\overline T$) and $W_1$ the Weyl group of $\overline R$. Then $W_1$ is $\sigma$-invariant, and so $\sigma$ acts on $N_W(W_1)/W_1$. Two elements $W_1w_1$ and $W_1w_2$ are $\sigma$-conjugate if there is $w\in W$ such that $W_1w_2=(W_1w)^\sigma(W_1w_1)(W_1w)^{-1}$.

\begin{lemma} \label{l:bijection}
The $\overline G_\sigma$-orbits of the set of $\sigma$-stable conjugates of $\overline R$ in $\overline G$ are in bijective
correspondence with $\sigma$-conjugacy classes in $N_W(W_1)/W_1$, with bijection inducing by $\overline R^g\mapsto W_1\pi(g^\sigma g^{-1})$.
\end{lemma}

\begin{proof}
See \cite[Proposition 3]{78Car}.
\end{proof}

\begin{lemma} \label{l:pperiod} Suppose that $\overline R^g$ is $\sigma$-stable. Then $\exp_p((\overline R^g)_\sigma)$ is equal to the minimal power of $p$ greater than the
maximal height of a root in $\Phi_1$.
\end{lemma}

\begin{proof}
This is a straightforward consequence of \cite[Proposition 0.5]{95Tes}.
\end{proof}

Denote by $mh(\Psi)$ the maximal height of a root in $\Psi$.  Subsystems of $E_7$ have type $A_l$, $D_l$, $E_6$, or $E_7$, and the corresponding
maximal heights are given in Table~\ref{tab:heights}.

\begin{table}[h]
\caption{Maximal height of a root in a root system $\Psi$}\label{tab:heights}
$\begin{array}{|c|c|c|c|c|}
\hline
\Psi&A_l&D_l&E_6&E_7\\
\hline
mh(\Psi)&l&2l-3&11&17\\
\hline
\end{array}$
\end{table}

\begin{lemma}\label{l:center} Suppose that $\overline R^g$ is $\sigma$-stable, $g\in N_{\overline G}(\overline T)$  and $\pi(g^\sigma g^{-1})=w$. Then $Z((\overline R^g)_\sigma)$ is conjugate in
$\overline G$ to $(Z(\overline R))_{\sigma\circ w}=Z(\overline R)\cap \overline T_{\sigma\circ w}$.
\end{lemma}

\begin{proof}
See \cite[Lemma 2.6]{22ButGre}.
\end{proof}

The next lemma is concerned with the structure of $N_W(W_1)/W_1$. Let $\Pi_1$ and $\Delta_1$ be the fundamental system and Dynkin diagram of $\Phi_1$ respectively. We write $N_W(\Pi_1)$ for the set of elements of $W$ such that $w(\Pi_1)=\Pi_1$ and $\Aut_W(\Delta_1)$ for the group of symmetries of $\Delta_1$ induced by transformations by elements of $W$. Also denote by $\Phi_2$  the set of roots orthogonal to every root of $\Phi_1$ and by $W_2$ the subgroup of $W$ generated by reflections in roots of $\Phi_2$. It is clear that $W_2\leq N_W(\Pi_1)$.

\begin{lemma}\label{l:N(W1)} $W_1\times W_2$ is a normal subgroup of $N_W(W_1)$, $N_W(W_1)=W_1N_W(\Pi_1)$ with $W_1\cap N_W(\Pi_1)=1$, $W_2=C_W(\Pi_1)$ and $N_W(\Pi_1)/W_2\simeq \Aut_W(\Delta_1)$.
\end{lemma}

\begin{proof}
See \cite[Proposition 28]{72CarW} and its proof.
\end{proof}

Given a root system $\Phi$ with a fundamental system $\Pi$, we denote by $w_0(\Pi)$ the unique
element of $W(\Phi)$ that maps every positive (with respect to $\Pi$) root of $\Phi$ to a negative one; this element has order 2 (see, for example, \cite[Proposition 2.2.6]{72Car}).
For a indecomposable system $\Phi$, we have $w_0(\Pi)=-1$ if and only if $\Phi\neq A_l$ ($l\geqslant 2$), $D_l$ ($l$ odd), $E_6$. In the following lemma, we regard $w_0(\Pi_1)$ as an element of $W$.

\begin{lemma}\label{l:-1} Suppose that $-1\in W$ and let $\theta$ be the natural homomorphism from $N_W(W_1)$ onto $N_W(\Pi_1)\simeq N_W(W_1)/W_1$. Then the following hold:

\begin{enumerate}
 \item $\theta(-1)$ is equal to $-w_0(\Pi_1)$ and acts as $-1$ on $\Phi_2$;
 \item $\theta(-1)\in W_2$ if and only if $w_0(\Pi_1)$ acts as $-1$ on $\Pi_1$, and if $\theta(-1)\not\in W_2$ then  $W_2\times \langle \theta(-1)\rangle\leqslant N_W(\Pi_1)$.
\end{enumerate}
\end{lemma}

\begin{proof}
 It is clear that $-w_0(\Pi_1)\in N_W(\Pi_1)$, and so (1) follows.  Now (2) follows from the facts that $W_2=C_W(\Pi_1)$ and $\langle \theta(-1)\rangle$ is a central subgroup of
 $N_W(\Pi_1)$.
\end{proof}

We conclude with the lemma on conjugacy classes of the Weyl groups $W(B_n)$ and $W(D_n)$. It is well known that every element $w$ of $W(B_n)$ is a product of positive and negative cycles and the signed cycle-type of $w$ is denoted by $[l_1,\dots l_j, \overline l_{j+1},\dots,\overline l_s]$, where $l_i$ and $\overline l_i$ stand for a positive cycle and a negative cycle of length $l_i$ respectively.
Also $-1\in W(B_n)$ and the signed cycles types $[l_1,\dots l_j, \overline l_{j+1},\dots,\overline l_s]$ and $[l'_1,\dots l'_j, \overline l'_{j+1},\dots,\overline l'_s]$ of $w$ and $-w$ respectively satisfy $l_i'=l_i$, $\overline l_i'=\overline l_i$ if $l_i$ is even and $l_i'=\overline l_i$, $\overline l_i'=l_i$ if $l_i$ is odd.

\begin{lemma}\label{l:dn}  An element of $W(B_n)$ lies in $W(D_n)$ if and only if
it has an even number of negative terms in its signed cycle-type.
Two elements of the Weyl group $W(D_n)$ are conjugate if and only if they have the same signed cycle-type, except that if all the cycles are even and positive there are two conjugacy classes. If $n$ is even and all the cycles in the signed cycle-type of $w$ are even and positive then $w$ and $-w$ are conjugate if and only $n/2$ is even.
\end{lemma}

\begin{proof}
See \cite[Proposition 25]{72CarW} and its proof.
\end{proof}

\section{Proof of Theorem \ref{t:E7}}\label{s:diag}

\subsection{Notation}

\begin{figure}[h]
\includegraphics{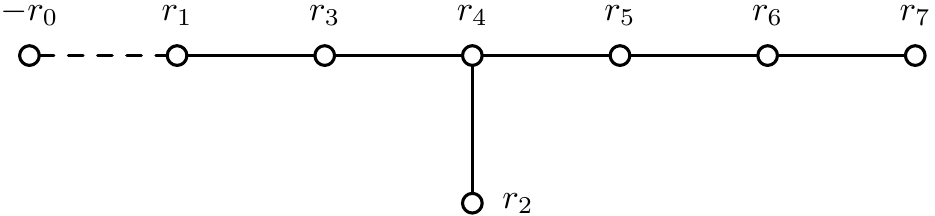}
\caption{Extended Dynkin diagram of $E_7$}\label{fig:E7}
\end{figure}

Let $p$ be an odd prime and $\overline F$   the algebraic closure of the field of  order $p$. Let $\overline L$ be the simply connected  simple algebraic group of type $E_7$ over $\overline F$, $\overline Z=Z(\overline L)$, and $\widetilde L=\overline L/\overline Z$.   We fix a maximal torus $\overline T$ of~$\overline L$, and let $\Phi$ be the root system of $\overline L$ with respect to~$\overline T$ and
$\Pi=\{r_1,\dots,r_7\}$ a fundamental subsystem of $\Phi$. Our numbering of fundamental roots is given in Fig. \ref{fig:E7},
where $r_0$ denotes the root of maximal height.

Let $x_r(s)$, $h_r(t)$ with  $r\in\Phi$, $s\in\overline F$, $t\in
\overline F^\times$ be Chevalley generators of  $\overline L$.  If $r=r_i$, then $h_i(t)=h_{r_i}(t)$.

For a power $q$ of $p$, let $\sigma$ be the endomorphism of $\overline L$ that maps $x_r(s)$ to $x_r(s^q)$. Denote the endomorphism of $\widetilde L$ induced by $\sigma$ by the same symbol. Then $Z(\overline L_\sigma)=\overline Z$ and $\overline L_\sigma/\overline Z$ is isomorphic to the simple group $E_7(q)$, so we identify $L$ with $\overline L_\sigma/\overline Z$. Also
by \cite[Lemma 2.5.8(a)]{98GorLySol}, it follows that $\widetilde L_\sigma$
is isomorphic to the group of inner-diagonal automorphisms of $E_7(q)$, so we put $G=\widetilde L_\sigma$.


As we explained in Section \ref{s:pre}, the spectrum of $G$ can be found as the set of all divisors of the numbers $\eta(R)$,
where $R$ runs over reductive subgroups of $G$ of maximal rank. By Lemma \ref{l:bijection}, we may assume that $R=(\widetilde R^g)_\sigma$, where $\widetilde R$ is a
$\sigma$-stable reductive subgroup of $\widetilde L$ containing $\widetilde T=\overline T/\overline Z$, $g\in N_{\widetilde L}(\widetilde T)$, $w=\pi(g^\sigma g^{-1})$ normalizers the Weyl group $W_1$ of $\widetilde R$,
and $W_1w$ runs over representatives of $\sigma$-conjugacy classes of $ N_W(W_1)/W_1$. Observe that $\sigma$ acts trivially on $W$, so the $\sigma$-conjugacy is the ordinary conjugacy.

By Lemmas \ref{l:pperiod} and \ref{l:center}, for such $R$ we have
\begin{equation}\label{e:eta} \eta(R)=n_p(mh(\Phi_1))\cdot \exp(Z(\widetilde R)_{\sigma\circ w}).\end{equation}
Thus $\eta(R)$ depends only on $\Phi_1$ and $w$, and sometimes we write $\eta(\Phi_1,w)$ instead of $\eta(R)$.

\subsection{The reductive subgroups to consider and their centers}
It is clear that to find $\omega(G)$, it suffices to consider only whose
numbers $\eta(\Phi_1,w)$ that are maximal under divisibility. In Table \ref{tab:SubsysE7}, we list
equivalence classes of subsystems of $E_7$ divided into two subsets $\alpha(\Phi)$ and
$\beta(\Phi)$. The next lemma shows that maximal under divisibility numbers arise from $\eta(\Phi_1,w)$ with $\Phi_1\in\alpha(\Phi)$.

\begin{table}[h]
\caption{Subsystems of $E_7$}\label{tab:SubsysE7}
\begin{tabular}{|c|l|}\hline
$\alpha(\Phi)$& $\varnothing$, $A_1$, $A_1^2$,  $(A_1^3)'$, $(A_1^3)''$, $A_3$, $(A_1^4)'$, $(A_1^4)''$, $D_4$, $A_1^5$, $(A_5)'$, $(A_5)''$,\\
& $D_5$, $A_3^2$, $D_6$, $E_6$,  $A_7$, $E_7$\\ \hline
$\beta(\Phi)$ & $A_2$, $A_2\times A_1$,  $A_2\times A_1^2$, $A_2^2$, $(A_3\times A_1)'$, $(A_3\times A_1)''$,  $A_4$,\\
&$A_2\times A_1^3$, $A_2^2\times A_1$, $(A_3\times A_1^2)'$, $(A_3\times A_1^2)''$, $A_3\times A_2$, $A_4\times A_1$,\\
& $D_4\times A_1$, $A_1^6$, $A_2^3$, $A_3\times A_1^3$, $A_3\times A_2\times A_1$, $A_4\times A_2$, $(A_5\times A_1)'$,\\
& $(A_5\times A_1)''$, $A_6$, $D_4\times A_1^2$, $D_5\times A_1$, $A_1^7$, $A_3^2\times A_1$, $A_5\times A_2$,\\
& $D_4\times A_1^3$, $D_6\times A_1$\\ \hline
\end{tabular}

\end{table}

\begin{lemma}\label{l:ListOfSubsystems} Every element of $\omega(G)$ divides $\eta(\Phi_1,w)$ for some $\Phi_1\in\alpha(\Phi)$
and $w$ permuting irreducible components of $\Phi_1$ cyclically.
%
\end{lemma}

\begin{prf} By \cite[Table III]{84Der}, reductive subgroups of maximal rank whose root system is $A_1^6$ or $A_1^7$ are not centralizers of semisimple elements.

Now suppose that $\Phi_1\in\beta(\Phi)$ with $\Phi_1\neq A_1^6$, $A_1^7$ and $w\in N_W(W_1)$.

First assume that $\Phi_1=A_l^i$ for some $l\in\{2,4,6\}$ and $1\leqslant i\leqslant 3$. We claim that for some $w_1\in W_1$, there is a $ww_1$-invariant subsystem of $\Phi_1$ of type $A_{l-1}^i$.
Let $\Psi_1$ be an irreducible component of $\Phi_1$ and take $\Psi_2$ to be some subsystem  of $\Psi_1$ of type $A_{l-1}$. Since
$\Psi_2w^i$ is also a subsystem of $\Psi_1$ of type $A_{l-1}$, there is $w_1\in W(\Psi_1)$ such that $\Psi_2w^iw_1=\Psi_2$. If $1\leqslant j<i$, then $\Psi_2w^jw_1=\Psi_2w^j$ because $w_1$ acts trivially on irreducible components of $\Phi_1$ other than $\Psi_1$. So $\Psi_2(ww_1)^i=\Psi_2$ and the orbit of $\Psi_2$ under the action of $\langle ww_1\rangle$ is the required invariant subsystem $\Phi_1'$ of $\Phi_1$ of type $A_{l-1}^i$. It is clear that $ww_1\in N_W(W_1')$, where $W_1'$ is the Weyl group of $\Phi_1'$.

Now $\eta(\Phi_1,w)=\eta(\Phi_1,ww_1)$ by Lemma \ref{l:bijection}. Also $n_p(mh(A_{l-1}))=n_p(mh(A_l))$ since $p$ is odd and $l$ is even. Finally, if $\widetilde R$ and  $\widetilde R'$ is the reductive groups corresponding to $\Phi_1$ and $\Phi_1'$ respectively, then $Z(\widetilde R)\leqslant Z(\widetilde R')$ and so
$Z(\widetilde R)_{\sigma\circ ww_1}\leqslant Z(\widetilde R')_{\sigma\circ ww_1}$. Thus $\eta(\Phi_1,w)$ divides $\eta(\Phi_1',ww_1)$.

Now if $\Phi_1\neq A_l^i$ for some $l\in\{2,4,6\}$ and $1\leqslant i\leqslant 3$, or $w$ does not permute irreducible components of $\Phi_1$ cyclically, then there is a $w$-invariant subsystem $\Phi_1'\in\alpha(\Phi)$ of $\Phi_1$ such that $n_p(mh(\Phi_1))=n_p(\Phi_1'))$ and the same argument as in the previous paragraph shows that $\eta(\Phi_1,w)$ divides $\eta(\Phi_1',w)$.
\end{prf}

\begin{table}[h]
\caption{Conditions on $Z(\overline R)$}\label{tab:CondE7}
\begin{tabular}{|c|c|c|}\hline
$\Phi_1$&$\Pi_1$& $Z(\overline R)$\\ \hline
$A_1$& $-r_0$&$t_1=1$\\ \hline
$A_1^2$& $-r_0$, $r_{11}$ & $t_1=t_6=1$\\ \hline
$(A_1^3)'$& $-r_0$, $r_{10}$, $r_{11}$& $t_1=t_4=t_6=1$\\ \hline
$(A_1^3)''$& $-r_0$, $r_7$, $r_{11}$ & $t_1=t_6=1$, $t_7^2=1$\\ \hline
$(A_1^4)'$& $-r_0$, $r_7$, $r_{10}$, $r_{11}$& $t_1=t_4=t_6=1$, $t_7^2=1$\\ \hline
$(A_1^4)''$& $-r_0$, $r_3$, $r_{10}$, $r_{11}$& $t_1=t_4=t_6=1$, $t_3^2=1$\\ \hline
$A_1^5$& $-r_0$, $r_5$, $r_7$, $r_{10}$, $r_{11}$& $t_1=t_4=t_6=1$, $t_5^2=t_7^2=1$\\ \hline
$A_3$& $-r_0$, $r_1$, $r_3$ & $t_1=t_3=t_4=1$\\ \hline
$A_3^2$& $-r_0$, $r_1$, $r_3$, $r_5$, $r_6$, $r_7$& $t_1=t_3=t_4=1$, $t_5=t_7^3$, $t_6=t_7^2$,  $t_7^4=1$\\ \hline
$D_4$& $-r_0$, $r_1$, $r_3$, $r_9$ & $t_1=t_3=t_4=1$, $t_7=t_5$\\ \hline
$(A_5)'$& $-r_0$,  $r_1$, $r_2$, $r_3$, $r_4$ & $t_1=t_3=t_4=1$, $t_5=t_2^{-1}$, $t_2^2=1$\\ \hline
$(A_5)''$& $-r_0$, $r_1$, $r_3$, $r_4$, $r_5$ & $t_1=t_3=t_4=1$, $t_2=t_5^{-1}$, $t_6=t_5^2$\\ \hline
$D_5$& $-r_0$, $r_1$, $r_3$, $r_4$, $r_9$ & $t_1=t_3=t_4=1$, $t_5=t_2^{-1}$, $t_7=t_2^{-1}$\\ \hline
$A_7$& $-r_0$, $r_1$, $r_3$, $r_4$, $r_5$, $r_6$, $r_7$ & $t_1=t_3=t_4=1$, $t_2=t_7$, $t_5=t_7^3$, $t_6=t_7^2$, $t_7^4=1$\\ \hline
$D_6$& $-r_0$,  $r_1$, $r_2$, $r_3$, $r_4$, $r_5$& $t_1=t_3=t_4=t_6=1$, $t_5=t_2$, $t_2^2=1$\\ \hline
$E_6$& $-r_0$, $r_1$, $r_3$, $r_4$, $r_5$, $r_8$& $t_1=t_3=t_4=1$, $t_2=t_5^{-1}$, $t_6=t_5^{2}$, $t_7=t_5$\\ \hline
$E_7$& $\Pi$ & $t_1=t_3=t_4=t_6=1$, $t_2=t_5=t_7$, $t_2^2=1$\\ \hline
\end{tabular}

\end{table}
\
In Table \ref{tab:CondE7}, we list nonempty subsystems of $\Phi$ lying in $\alpha(\Phi)$.
The first and second columns give, respectively, the type of the subsystem and  the chosen fundamental system. Recall that $r_0$ is the root of maximal height, that is,
$r_0=(2,2,3,4,3,2,1)$, where $(a_1,\dots,a_7)$ stands for the root $a_1r_1+\dots+a_7r_7$.
Also $r_{8}=(0,1,0,1,1,1,0)$, $r_{9}=(0,1,1,2,2,1,0)$, $r_{10}=(0,-1,-1,-2,-1,0,0)$, and $r_{11}=(0,-1,-1,-2,-2,-2,-1)$.

The last column gives the conditions for $h\in\overline T$
to lie in $Z(\overline R)$, where $\overline R$ is the reductive subgroup of $\overline L$ generated by $\overline T$ and the root subgroups corresponding to  the roots in $\Phi_1$. Using the formula $$x_r(s)^{h_i(t)}=x_r(st^{\langle
r,r_i\rangle}),$$ where $\langle r,r_i\rangle=2(r,r_i)/(r_i,r_i)$, we see that $h=h_1(t_1)\dots h_k(t_7)$ centralizes $x_r(s)$ for all $r\in\Pi_1$ if and only if
$$t_1^{\langle r,r_1\rangle}\dots t_7^{\langle
r,r_7\rangle}=1\text { for all }r\in\Pi_1,$$
and some equivalents of these equations are listed in the last column. Also note that $\widetilde R=\overline R/\overline Z$ and $Z(\widetilde R)=Z(\overline R)/\overline Z$ since $x_r(s)^hx_r(s)^{-1}$ is a unipotent element, so it can lie in $\overline Z$ only if it is trivial.

It is clear that $Z(\overline R)$ includes the torus generated by $h_r(t)$ with $r\in\Phi_2$, and we denote this torus by $\overline T_2$.


\subsection{Direct calculation of $\exp(Z(\widetilde R)_{\sigma\circ w})$}
In this subsection, we explain how one can directly calculate the exponent of $Z(\widetilde R)_{\sigma\circ w}$ for a given $w\in N_W(W_1)$. Denote the preimages of $Z(\widetilde R)_{\sigma\circ w}$ and $\widetilde T_{\sigma\circ w}$ in $\overline L$ by $H_w$ and $T_w$ respectively.
It is clear that $H_w=Z(\overline R)\cap T_w$. If $r_iw=a_{i1}r_1+\dots+a_{i7}r_7$, then $$h_i(t)^{\sigma\circ w}=h_1(t^{q a_{i1}})\dots
h_k(t^{ q a_{i7}}).$$ So $T_w$ consists of the elements satisfying
\begin{equation}\label{e:h} h_1(t_1)\dots h_k(t_k)=h_1(t_1^{\varepsilon qa_{11}}\dots t_k^{\varepsilon qa_{k1}})\dots
h_k(t_1^{\varepsilon qa_{1k}}\dots t_k^{\varepsilon qa_{kk}})z^i\end{equation} for $i=0,1$.
By the last row of Table \ref{tab:CondE7}, we have  $z^i=h_1(t_0^{v_1})\dots h_7(t_0^{v_7})$, where $t_0^2=1$ and $v=(0,1,0,0,1,0,1)$.
Therefore, if
we define $M_w=(a_{ij})$, then $h_1(t_1)\dots h_7(t_7)\in T_w$ if and only if, in additive notation,
\begin{equation}\label{e:mat}(t_1,\dots,t_7,t_0)\begin{pmatrix}qM_w-E& \text{\large0}\\
v & 2\\
\end{pmatrix}=0.\end{equation}

Since we are interested only in elements of $Z(\overline R)\cap T_w$, we can reduce the number of equations. Suppose that $t_i=t^{(i)}c$ for elements of $Z(\overline R)$,
where $t^{(i)}=(t_1,\dots,t_{i-1},t_{i+1},\dots,t_7)$ and $c$ is an integer column vector of size $6$. Since $Z(\overline R)$ is
$w$-invariant, it follows that $((t_1,\dots,t_7)M_w)_i=((t_1,\dots,t_7)M_w)^{(i)}c$  and also $v_i=v^{(i)}c$. So the equation $t_i= q((t_1,\dots,t_7)M_w)_i+v_i$ is a corollary of $t_i=t^{(i)}c$ and $(t)^{(i)}=
q((t_1,\dots,t_k)M_w)^{(i)}+v^{(i)}$. Thus we can remove the $i$th column of the matrix in \eqref{e:mat}. Also after adding the $i$th row multiplied by $c_j$ to the $j$th row for all $j\neq i$, we can remove the $i$th row. We  repeat this procedure for all relations of the form $t_i=t^{(i)}c$ and denote the matrices obtained from $qM_w-E$  and $v$ by $X$ and $y$ respectively (observe that $y$ does not depend on $w$).

Sometimes we resolve the reduced system by hand, but in most cases we find, with help of \textsc{Magma} \cite{MagmaBook}, invertible integer matrices $A$ and $B$ such that $AXB$ is diagonal and $yB$ has a single non-zero entry modulo $2$. These matrices $A$, $B$, $AXB$,
as well as $yB$ modulo $2$, are given in Tables \ref{tab:7a1}-\ref{tab:7d5}.
The following example shows how the exponent of $H_w/\overline Z$ can be
readily seen from $A$, $B$, $AXB$ and $yB$.

Let $\Phi_1=A_1$ and $\Pi_1=\{-r_0\}$. Let
$$M_w=\begin{pmatrix}
 1 & 0 & 1 & 1 & 1 & 1 & 1\\
 0 & 1 & 0 & 1 & 1 & 1 & 1\\
 0 & 0 &-1 &-1 &-1 &-1 &-1\\
 0 & 0 & 1 & 0 & 0 & 0 & 0\\
 0 & 0 & 0 & 1 & 0 & 0 & 0\\
 0 & 0 & 0 & 0 & 1 & 0 & 0\\
 0 & 0 & 0 & 0 & 0 & 1 & 0\\
\end{pmatrix}.$$

The relation of $Z(\overline R)$ is $t_1=1$.  Thus removing the first row and column of $qM_w-E$ and the first entry of $v$, we result in
$$X=\begin{pmatrix}
q - 1&0&q&q&q&q\\
0&-q-1&-q&-q&-q&-q\\
0&q&-1&0&0&0\\
0&0&q&-1&0&0\\
0&0&0&q&-1&0\\
0&0&0&0&q&-1\\
\end{pmatrix}\text{ and }y=(1,0,0,1,0,1).$$

Let $A$ and $B$ be as in the first part of Table \ref{tab:7a1}
and $(s_2,\dots,s_7)=(t_2,\dots, t_7)A^{-1}$. Then $AXB=D=\diag(-2,(q^6-1)/2,-1,-1,1,1)$ and
$$(s_2,\dots,s_7,t_0)\begin{pmatrix}D& \text{\large0}\\
yB & 2\end{pmatrix}=0.$$

If $q\equiv 3\pmod 4$, then $yB\equiv (1,0,0,0,0,0)\pmod 2$, so $s_2^2t_0=1$, $s_3^{(q^6-1)/2}=1$ and hence $$H_w\simeq \mathbb{Z}_4\times \mathbb{Z}_{(q^6-1)/2}.$$
Also $y\equiv (1,1,0,0,0,0)A\pmod 2$, which means that $z$ lies in neither of the direct factors. It follows that the exponent of $H_w/\overline Z$ is equal to $[4,(q^6-1)/2]=(q^6-1)/2$.

If $q\equiv 1\pmod 4$, then $yB\equiv (1,1,0,0,0,0)\pmod 2$, so $s_2^2t_0=1$, $s_3^{(q^6-1)/2}t_0=1$. It follows that $H_w$ is generated by cyclic groups of orders $4$ and $q^6-1$. On the other hand, since  $y\equiv (0,1,0,0,0,0)A\pmod 2$, the cyclic group of order $q^6-1$ contains $z$, and hence the exponent of $H_w/\overline Z$ is equal to $(q^6-1)/2$, as in the previous case.

\normalsize

\subsection{The set $\omega(G)$} To prove that the elements of $\omega(G)$ are precisely
the divisors of elements of $\nu(q)$, it suffices to show that
for every $\Phi_1\in \alpha(\Phi)$ and every conjugacy class  $[w]$ of $N_W(\Pi_1)$, the number $\eta(\Phi_1,w)$ divides some element of $\nu(q)$ and, conversely,
every element of $\nu(q)$ arises as $\eta(\Phi_1,w)$ for some $\Phi_1$ and $w$.

\emph{Remark 1}. If  $g\overline Z\in Z(\widetilde R)_{\sigma\circ w}$, then $g^2\in \overline R_{\sigma\circ w}$, so the exponent of $Z(\widetilde R)_{\sigma\circ w}$ is at most twice that of  $Z(\overline R)_{\sigma\circ w}$ (or even that of $Z(\overline R)_{\sigma\circ w}/\overline Z$).

\emph{Remark 2}. Observe that $-1\in W$ and recall from Section \ref{s:pre}, that $\theta$ is the natural homomorphism from $N_W(W_1)$ onto $N_W(\Pi_1)\simeq N_W(W_1)/W_1$.
If $w\in N_W(\Pi_1)$, then $\theta(-1) w\in N_W(\Pi_1)$, and the exponent of $Z(\widetilde R)_{\sigma\circ(\theta(-1)w)}$ can be obtained from that of $Z(\widetilde R)_{\sigma\circ w}$ by replacing $q$ with $-q$, so if $w$ and $\theta(-1)w$ are not conjugate in $N_W(\Pi_1)$, then it sufficient to consider the conjugacy class of only one of them.

(i) Let $\Phi_1=\varnothing$, that is, $\widetilde R$ is a maximal torus of $\widetilde L$. Since $G$ and $\overline L_\sigma$ are in duality (in the sense of \cite[Proposition 4.3.1]{85Car}), they have isomorphic maximal tori by \cite[Proposition 4.4.1]{85Car}. The structure of maximal tori of $\overline L_\sigma$, the universal version of $E_7(q)$, is found in \cite{91DerFak}. Maximal (under divisibility) orders of semisimple elements derived from this information are given in \cite[Theorem 1]{16But.t}. These numbers are reproduced in Item (1) of Theorem \ref{t:E7}.


(ii) Let $\Phi_1=A_1$ and $\Pi_1=\{-r_0\}$. Then $\Phi_2$ is of type $D_6$ with fundamental roots $r_i$, $2\leqslant i\leqslant 7$ and $N_W(\Pi_1)=W_2$ by \cite[p. 195]{83Der}. Also by Table~\ref{tab:CondE7}, we see that $Z(\overline R)=\overline T_2$ is a maximal torus of the subsystem subgroup $\overline Q$ of $\overline L$ corresponding to $\Phi_2$. So the groups $Z(\overline R)_{\sigma\circ w}$ are isomorphic to maximal tori of $\overline Q_\sigma$. Since $\overline Q_\sigma\simeq Spin_{12}^+(q)$, the exponents of these tori can be deduced from \cite{15Zav.t}. Doubling these exponents, we see that
$\exp(Z(\widetilde R)_{\sigma\circ w})$ divides one of the numbers $q^6-1$, $2(q^5\pm 1)$, $(q^4+1)(q^2\pm 1)$, $(q^2+1)(q^3\pm 1)(q\mp1)$, and $2(q^4-1)$.
We claim that the exponents $q^6-1$, $2(q^5\pm 1)$ and $2(q^4-1)$ do not arise, while the exponents $(q^4+1)(q^2\pm 1)$ and $(q^4-1)(q^2\pm q+1)$ do.

By Lemma \ref{l:-1}, we see that $\theta(-1)\in W_2$ and
it acts as $-1$ on $\Phi_2$, so we can use Lemma \ref{l:dn}
to decide whether $w$ and $\theta(-1)w$ are conjugate.

It follows from \cite[Theorem 1]{15Zav.t} that a torus of exponent $q^6-1$, $2(q^5-1)$, or $2(q^5+1)$ would correspond to $w\in W_2$ of type
$[6]$, $[5,1]$, or $[\overline 5,\overline 1]$ respectively (if we identify $W_2$ with $W(D_6)$ and use the notation introduced before Lemma \ref{l:dn}). Similarly, a torus of exponent $2(q^4-1)$ would correspond to an element of type $[4,2]$, $[4,1,1]$, or $[4,\overline 1,\overline 1]$. By the preceding paragraph, it is sufficient to carry out calculations only for elements of types $[6]$, $[4,2]$, $[4,1,1]$ and $[5,1]$, and the  corresponding results are presented in the first four parts of Table \ref{tab:7a1}. The remaining parts of Table  \ref{tab:7a1} show that the exponents $(q^4+1)(q^2\pm 1)$ and $(q^2+1)(q^3-1)(q+1)$ do arise. Thus we have all numbers in Item (2) except for $p(q^4-q^2+1)$.

\tiny

\begin{table}

\caption{$\Phi_1=A_1$, $\Pi_1=\{-r_0\}$, $t_1=1$, $y=(1,0,0,1,0,1)$}
\label{tab:7a1}

$\begin{array}{|rclrcl|}
\hline
M_w&=&\begin{pmatrix}
 1 & 0 & 1 & 1 & 1 & 1 & 1\\
 0 & 1 & 0 & 1 & 1 & 1 & 1\\
 0 & 0 &-1 &-1 &-1 &-1 &-1\\
 0 & 0 & 1 & 0 & 0 & 0 & 0\\
 0 & 0 & 0 & 1 & 0 & 0 & 0\\
 0 & 0 & 0 & 0 & 1 & 0 & 0\\
 0 & 0 & 0 & 0 & 0 & 1 & 0\\
\end{pmatrix}, & X&=&\begin{pmatrix}
q - 1&0&q&q&q&q\\
0&-q-1&-q&-q&-q&-q\\
0&q&-1&0&0&0\\
0&0&q&-1&0&0\\
0&0&0&q&-1&0\\
0&0&0&0&q&-1\\
\end{pmatrix},\\
A&=&\begin{pmatrix}
1 &  1 &  0 &  0 &  2 &  2\\
\frac{q^6-1}{2(q-1)} & \frac{q^6-2q+1}{2(q-1)} & \frac{q^5-q}{q-1}&  \frac{q^4-q}{q-1}& q(q+1)& q\\
0 &  0 &  1 &  0 &  0 &  0\\
0 &  0 &  0 &  1 &  0 &  0\\
0 &  0 &  0 &  0 &  1 &  0\\
0 &  0 &  0 &  0 &  0 &  1\\
\end{pmatrix},&
B&=&\begin{pmatrix}
  1 &  \frac{q+1}{2} &  0 &  0 &  1 &  1\\
1 &  \frac{q- 1}{2} &  0 &  0 &  1 &  1\\
q &  \frac{q^2-q}{2} &  1 &  0 &  q &  q\\
q^2 &  \frac{q^3-q^2}{2} &  q &  1 &  q^2 &  q^2\\
q^3 &  \frac{q^4 - q^3}{2} &  q^2  & q  & q^3 - 1  & q^3\\
q^4 &  \frac{q^5 - q^4}{2} &  q^3 &  q^2 &  q^4 - q &  q^4 - 1\\
\end{pmatrix},\\
q&\equiv& 1\mypmod 4: y\equiv A[2]\mypmod2,&yB&\equiv& (1,1,0,0,0,0)\mypmod2\\
q&\equiv& 3\mypmod 4: y\equiv A[1]+A[2]\mypmod2,&yB&\equiv& (1,0,0,0,0,0)\mypmod2\\
AXB&=&\diag(-2,(q^6-1)/2,-1,-1,1,1),&\exp&=&[4, (q^6-1)/2]=(q^6-1)/2\\

\hline
M_w&=&\begin{pmatrix}
 1 & 0 & 1 & 0 & 0 & 0 & 0\\
 0 & 1 & 0 & 0 & 0 & 0 & 0\\
 0 & 0 &-1 & 0 & 0 & 0 & 0\\
 0 & 0 & 1 & 1 & 1 & 1 & 1\\
 0 & 0 & 0 & 0 &-1 &-1 &-1\\
 0 & 0 & 0 & 0 & 1 & 0 & 0\\
 0 & 0 & 0 & 0 & 0 & 1 & 0\\
\end{pmatrix},&
X&=&\begin{pmatrix}
 q - 1   &   0  &    0  &    0  &    0 &     0\\
     0 &-q - 1  &    0 &     0  &    0&      0\\
     0  &    q & q - 1  &    q &     q &     q\\
     0  &    0 &     0 &-q - 1 &    -q &    -q\\
     0   &   0  &    0 &     q &    -1 &     0\\
     0  &    0  &    0  &    0  &    q&     -1\\
\end{pmatrix},\\
A&=&\begin{pmatrix}
1 &  0 &  0 &  0  & -q + 1 &  -q + 1\\
0 &  1  & 0 &  0 &  0 &  0\\
0  & q^3 + q &  \frac{q^4-1}{q-1} &  \frac{q^4-q}{q-1} &  q^2 + q &  q\\
0 &  0  & 1 &  1 &  0 &  0\\
0 &  0 &  0 &  0  & 1  & 0\\
0 &  0 &  0 &  0 &  0 &  1\\
\end{pmatrix},&
B&=&\begin{pmatrix}
1     &    0  &       0    &     0   &     -1    &    -1\\
0     &    1   &      0    &     0   &      0    &     0\\
0     &    0  &       1   &      0   &      0   &      0\\
0    &     q  &   q - 1   &      1   &      0    &     0\\
0    &   q^2 &  q^2 - q  &       q    &     1   &      0\\
0    &   q^3& q^3 - q^2   &    q^2    &     q    &     1\\
\end{pmatrix},\\
y&\equiv &A[1]+A[3]\mypmod2,& yB&\equiv& (1,0,0,0,0,0)\mypmod2,\\
AXB&=&\diag(q-1,-q-1,q^4-1,1,1,1),&\exp&=&[2(q-1),q+1,q^4-1]=q^4-1\\
\hline
M_w&=&\begin{pmatrix}
 1 & 0 & 0 & 0 & 0&  0 & 0\\
 0 & 1 & 0 & 0&  0&  0 & 0\\
 0 & 0 & 1 & 0 & 0&  0&  0\\
 0 & 0 & 0 & 1 & 1 & 1 & 1\\
 0 & 0 & 0&  0& -1& -1& -1\\
 0 & 0&  0&  0&  1&  0&  0\\
 0 & 0&  0&  0&  0&  1&  0\\
\end{pmatrix},&
X&=&\begin{pmatrix}
 q - 1  &    0 &     0 &     0  &    0 &     0\\
     0  &q - 1 &     0 &     0  &    0 &     0\\
     0  &    0 & q - 1 &     q  &    q &     q\\
     0  &    0 &      0& -q - 1 &    -q&     -q\\
     0  &    0 &     0 &     q &    -1 &     0\\
     0  &    0 &     0 &     0 &     q &    -1\\
\end{pmatrix},\\
A&=&\begin{pmatrix}
0 &  0  & 0 &  0  & 1 &  0\\
0  & 0  & 0 &  1 &  0 &  0\\
0  & 0 &  1 &  0 &  0  & 0\\
0  & -1 &  0 &  q - 1  & 0 &  0\\
-1 &  0  & 0 &  -q + 1 &  q - 1&   0\\
0  & 0 &  \frac{q^4-1}{q-1} &  \frac{q^3-1}{q-1} &  q + 1 &  1\\
\end{pmatrix},&
B&=&\begin{pmatrix}
1   &       -1    &       0      &     0    &      -1    &       0\\
   0    &       1    &       0    &      -1    &       0    &       0\\
       0   &       0     &    -1  &        0  &        0  &        q\\
        0     &    -1      &   -q     &     0     &     0   & q^2 - q\\
       -1  &       -q    &   -q^2    &      0    &      0  &q^3 - q^2\\
        1   &   q + 1&q^2 + q + 1  &        0    &      0   &-q^3 + 1\\
\end{pmatrix},\\
y&\equiv& A[5]+A[6]\mypmod 2,& yB&\equiv& (0,0,0,0,1, 0)\mypmod2\\
AXB&=&\diag(1,1,1,q-1, q-1,q^4-1),& \exp&=&[2(q-1),q^4-1]=q^4-1\\
\hline
\end{array}$
\end{table}

\addtocounter{table}{-1}

\begin{table}

\caption{continued}

$\begin{array}{|rclrcl|}
\hline
M_w&=&\begin{pmatrix}
 1 & 0 & 0 & 0 & 0 & 0 & 0\\
 0 & 1 & 0 & 1 & 1 & 1 & 1\\
 0 & 0 & 1 & 1 & 1 & 1 & 1\\
 0 & 0 & 0 &-1 &-1 &-1 &-1\\
 0 & 0 & 0 & 1 & 0 & 0 & 0\\
 0 & 0 & 0 & 0 & 1 & 0 & 0\\
 0 & 0 & 0 & 0 & 0 & 1 & 0\\
\end{pmatrix},&
X&=&\begin{pmatrix}
 q - 1  &    0 &     q  &    q  &    q   &   q\\
     0  &q - 1   &   q   &   q   &   q   &   q\\
     0 &     0& -q - 1   &  -q &    -q   &  -q\\
     0 &     0 &     q  &   -1   &   0   &   0\\
     0  &    0  &    0  &    q &    -1    &  0\\
     0  &    0  &    0   &   0  &    q    & -1\\
\end{pmatrix},\\
A&=&\multicolumn{4}{l|}{\begin{pmatrix}
1 &  0 &  0 &  0 &  0 &  0\\
0 &  0 &  0 &  0 &  1 &  0\\
0 &  0 &  1 &  0 &  0 &  0\\
0 &  0 &  0 &  1 &  0 &  0\\
1 &  -1 &  0 &  0 &  q - 1&   0\\
\frac{q^5-1}{q-1} &  0 &  \frac{q^4-1}{q-1} &  q^5 + q^2 + q &  q^5 + q &  1\\
\end{pmatrix},}\\
B&=&\multicolumn{4}{l|}{\begin{pmatrix}
 -1  & q &  0 &  -q &  0 &  q\\
-1 &  q - 1 &  0  & -q  & -1 &  q\\
-q &  q^2 - q  & -1 &  -q^2 + q &  0 &  q^2 - q\\
-q^2 &  q^3 - q^2 &  -q &  -q^3 + q^2 - 1 &  0 &  q^3 - q^2\\
-q^3 &  q^4 - q^3 - 1 &  -q^2  & -q^4 + q^3 - q &  0 &  q^4 - q^3\\
\frac{q^4-1}{q-1} &  -q^4 &  q^2 + q + 1 &  q^4 + q &  0 &  -q^4 + 1\\
\end{pmatrix},}\\
y&\equiv&A[6]\mypmod 2,& yB&\equiv& (0,0,0,0,0, 1)\mypmod2,\\
AXB&=&\diag(1,-1,1,1, q-1,q^5-1),& \exp&=&q^5-1\\
\hline
M_w&=&\begin{pmatrix}
 1 & 0 & 1 & 1 & 1 & 1 & 1\\
 0 & 0& -1 & 0 & 0 & 0 & 0\\
 0 & 1 & 0 & 0 & 0 & 0 & 0\\
 0 &-1 & 0 &-1 &-1 &-1 &-1\\
 0 & 1 & 1 & 2 & 1 & 1 & 1\\
 0 & 0 & 0 & 0 & 1 & 0 & 0\\
 0 & 0 & 0 & 0 & 0 & 1 & 0\\
\end{pmatrix},&
X&=&\begin{pmatrix}
    -1  &   -q  &    0  &    0   &   0  &    0\\
     q  &   -1  &    0  &    0   &   0  &    0\\
    -q  &    0 &-q - 1  &   -q  &   -q  &   -q\\
     q  &    q &   2q & q - 1  &    q   &   q\\
     0  &   0   &   0 &     q &    -1 &     0\\
     0  &    0  &    0&      0&      q&     -1\\
\end{pmatrix},\\
A^\top&=&\multicolumn{4}{l|}{\begin{pmatrix}
1 & 0 &  0 &  \frac{q^6+q^4+q^2+2q+1}{2} &  0 &  1\\
0 &  1&   q&   \frac{-q^5+q^4+q+ 1}{2}&   0 &  1\\
0 &  1&   q + 1 & -q^5 + q &  0 &  1\\
0 & 0&   1&   \frac{-q^5-q^4+q+1}{2} &  0 &  0\\
0 &  -q^2-q  & -q^3-q^2 &  \frac{q^7+q^6-q^5+q^4-q^3+q^2+q+1}{2} &  1 &  -q^2-q+1\\
0 &  -q  & -q^2 &   \frac{2q^6-q^5 + q^4+q + 1}{2} &  0&-q + 1\\
\end{pmatrix},}\\
B&=&\multicolumn{4}{l|}{\begin{pmatrix}
1&-q&-q^4 - q^3 + q^2 + q&-\frac{q^5 + q}{2}&0&\frac{q^5+q}{2}\\
0&   1 &  q^3 + q^2 - q - 1 &\frac{q^4 + 1}{2}& 0 & -\frac{q^4+1}{2}\\
0&   0&   -q^3 - 2q^2 - q + 1 &  -\frac{q^4+q^3+q^2+1}{2} &  0  & \frac{q^4 + q^3 + q^2 + 1}{2}\\
0 &  0 &  q + 1 &  \frac{q^2 + 1}{2} &  0 &  -\frac{q^2 + 1}{2}\\
0 &  0 &  q^2 + q &  \frac{q^3 + q}{2} &  1 &  -\frac{q^3 + q}{2}\\
-1&   -1&   q^3 + q^2  & \frac{q^4 + q^2}{2}&   q - 1 &  -\frac{q^4 + q^2-2}{2}\\
\end{pmatrix},}\\
y&\equiv& A[3]+A[4]\mypmod2,& yB&\equiv& (0,0,0,1,0,0)\mypmod2\\
AXB&=&\diag(1,1,2,(q^4+1)(q^2+1)/2,1,1),& \exp&=&(q^4+1)(q^2+1)\\
\hline
\end{array}$
\end{table}

\addtocounter{table}{-1}

\begin{table}

\caption{continued}

$\begin{array}{|rclrcl|}
\hline
M_w&=&\begin{pmatrix}
 1 & 1 & 1 & 1 & 1 & 1 & 1\\
 0 & 0& -1 & 0 & 0 & 0 & 0\\
 0 &-1&  0&  0&  0&  0&  0\\
 0 & 0&  0& -1& -1& -1& -1\\
 0 & 1&  1&  2&  1&  1&  1\\
 0 & 0&  0&  0&  1&  0&  0\\
 0 & 0&  0&  0&  0&  1&  0\\
\end{pmatrix},&
X&=&\begin{pmatrix}
    -1 &    -q&      0&      0&      0&      0\\
    -q &    -1&      0&      0&      0&      0\\
     0 &     0& -q - 1&     -q&     -q&     -q\\
     q &     q&    2q&  q - 1&      q&      q\\
     0 &     0&      0&      q&     -1&      0\\
     0 &     0&      0&      0&      q&     -1\\
\end{pmatrix},\\
A^\top&=&\multicolumn{4}{l|}{\begin{pmatrix}
1&q^3 + q^2 + q&0&0&-q^3 - q^2 - q&\frac{-q^8+q^5+q^4+q}{2}\\
0& 0&0&0&0&-1\\
0 &  q^3 + q^2 + q + 1 &  1 &  0&   -q^3 - q^2 - q - 1&   \frac{-q^8+2q^4-1}{2}\\
0 &  q^2 + q + 1&   0 &  1&   -q^2 - q&   \frac{-q^7+q^4+q^3-1}{2}\\
0 &  q+1 &  0 &  0 &  -q&   \frac{-q^6+q^4+q^2-1}{2}\\
0&   1  & 0&  0& -1& \frac{-q^5 +q^4 + q - 1}{2}\\
\end{pmatrix},}\\
B&=&\multicolumn{4}{l|}{\begin{pmatrix}
-1 &  \frac{q^4-q^3+q^2-q}{2}&   0 &  0 &  0 &  q^5 + q\\
0  & \frac{-q^3+q^2-q+1}{2}& 0&0&0&-q^4 - 1\\
0  & \frac{-q^4 +q^3 -q^2 + 2q}{2} &  -1 &  0 &  0 &  -q^5 + q^2\\
-q &  \frac{q^2 - q}{2} & -q &  -1 &  0 &  q^3 - q\\
-q^2&  \frac{q^3-q^2 -2}{2}&   -q^2 &  -q + 1 &  -1 &  q^4 - q^2\\
q^2 + q &  \frac{q^4 - q^3}{2}  & q^2 + q + 1  & q  & 1 &  q^5 - q^3\\
\end{pmatrix},}\\
y&\equiv& A[5]+A[6]\mypmod 2, & yB&\equiv& (0,0,0,0,1,0)\mypmod2\\
AXB&=&\diag(1,1,1,1,1,(q^4+1)(q^2-1)),& \exp&=&(q^4+1)(q^2-1)\\
\hline
M_w&=&\begin{pmatrix}
  1 & 0 & 1 & 1 & 0&  0&  0 \\

  0 &-1 & 0& -1&  0&  0&  0 \\

  0 & 0& -1& -1&  0&  0&  0 \\

  0 & 1&  1&  1&  0&  0&  0 \\

  0 & 0&  0&  1&  1&  1&  1 \\

  0 & 0&  0&  0&  0& -1& -1 \\

  0 & 0&  0&  0&  0&  1&  0 \\\end{pmatrix},&
X&=&\begin{pmatrix}
 -q-1 & 0& -q&  0&  0&  0 \\

  0& -q-1& -q&  0&  0&  0 \\

  q&  q&  q-1&  0&  0&  0 \\

  0&  0&  q&  q-1&  q&  q \\

  0&  0&  0&  0& -q-1& -q \\

  0&  0&  0&  0&  q&  -1 \\
\end{pmatrix},\\
A^\top&=&\multicolumn{4}{l|}{\begin{pmatrix}
 1  & 0&    \frac{-q^7 - q^6 + 2q^3 + 2q^2 + 3q + 1}{2}&   q^2&   0&   0 \\

 0  & 1&  \frac{-q^6 - q^5 - q^4 +

    3q^2 + 3q + 3}{2}&   0&   0&   0 \\

1&0    &  \frac{-q^7 - 2q^6 + q^4 + 3q^3 + 4q^2 + 4q +

    3}{2} & q^2+q
    &0&0  \\

 0 &  0 &  \frac{-q^7 - 2q^6 - q^5 + q^4 + 4q^3 + 5q^2 + 4q + 2}{2}
&  q^2+q+1 &  1&   0 \\

 0  & 0 &  \frac{-q^7 - 2q^6 + q^4 + 3q^3 + 4q^2 + 2q + 1}{2} &  q^2+q &  1 &  0 \\

 0  & 0&   \frac{-2q^6 - q^5 +q^3 + 3q^2 + 2q + 1}{2} &  q &  0 &  1 \\
\end{pmatrix},}\\
B^\top&=&\multicolumn{4}{l|}{\begin{pmatrix}\frac{q^5+q^4+q^3-q^2-q-3}{2}&\frac{q^4-q}{2}&\frac{-q^4-q^3+q+1}{2}&\frac{q^3+q^2}{2}&\frac{-q^5+q}{2}&\frac{-q^6+q^2}{2}\\
\frac{q^5+q^4+q^3-q^2-3q-3}{2}&\frac{q^4-q-2}{2}&\frac{-q^4-q^3+q+3}{2}&\frac{q^3+q^2}{2}&\frac{-q^5+3q}{2}&\frac{-q^6+3q^2}{2}\\
\frac{q^5+q^4+q^3-q^2-q-1}{2}&\frac{q^4-q}{2}&\frac{-q^4-q^3+q+1}{2}&\frac{q^3+q^2}{2}&\frac{-q^5+q}&\frac{-q^6+q^2}{2}\\
-q^4+q^2+2q+1&-q^3+q^2+q&q^3-2q-1&-q^2+2&q^4-q^3-q^2+q-2&q^5-q^4-q^3+q^2-2q\\
\frac{q^5+q^4+q^3-q^2-q-1}{2}&\frac{q^4-q}{2}&\frac{-q^4-q^3+q+1}{2}&\frac{q^3+q^2}{2}&\frac{-q^5+q-2}{2}&\frac{-q^6+q^2+2q}{2}\\
\frac{q^5+q^4+q^3-q^2-q-1}{2}&\frac{q^4-q}{2}&\frac{-q^4-q^3+q+1}{2}&\frac{q^3+q^2}{2}&\frac{-q^5+q}{2}&\frac{-q^6+q^2+2}{2}\\
\end{pmatrix},}\\
y&\equiv& A[4]\mypmod2,& yB&\equiv&(0,0,1,0,0,0)\mypmod2,\\
AXB&=&-\diag(1,1,(q^4-1)(q^2+q+1)/2,-2,1,1),& \exp&=&(q^4-1)(q^2+q+1)=(q^3-1)(q^2+1)(q+1)\\
\hline
\end{array}$
\end{table}

\normalsize

(iii) Suppose that $\Phi_1=A_1^k$ and $\Phi_1\neq (A_1^3)''$.  Using \cite[Corollary 1]{16But.t}, where the numbers $\exp(Z(\overline R)_{\sigma\circ w}/\overline Z)$ were calculated, we see that every $\eta(\Phi_1,w)$ divides one of the numbers $p(q^4-1)$, $2p(q^2-1)(q^2\pm q+1)$, and $p(q^4+1)(q\pm 1)$, and so give nothing new.

(iv) Let $\Phi_1=(A_1^3)''$ and $\Pi_1=\{-r_0,r_7,r_{11}\}$. Then $\Phi_2$ has type $D_4$ with fundamental roots $r_2$, $r_3$, $r_4$, and $r_5$, and $Z(\overline R)$ consists of the elements satisfying $t_1=t_6=1$, $t_7^2=1$. It follows that $Z(\overline R)=\overline T_{2}\times \overline Z$. Since $w$ stabilizes $\Phi_2$, we conclude that $H_w=(\overline T_2)_{\sigma\circ w}\times \overline Z$, and so $H_w/\overline Z\simeq  (\overline T_{2})_{\sigma\circ w}$.

By \cite[p. 197]{83Der}, the group $N_W(\Pi_1)$ induces the full automorphism group on the Dynkin diagram of $\Phi_2$, so $(\overline T_{2})_{\sigma\circ w}$, where $w$ runs over conjugacy classes
of $N_W(\Pi_1)$, are exactly the maximal tori of $Spin_8^+(q)$, $Spin_8^-(q)$ and $^3D_4(q)$. The maximal exponents of these tori are $q^4-q^2+1$, $q^4+1$, $(q^4-1)/2$ and $(q^2-1)(q^2\pm q+1)$ (see \cite[Table 1]{15Zav.t} for the spin groups and \cite[Table 1.1]{87DerMich} for $^3D_4(q)$). The first number gives $p(q^4-q^2+1)$ in Item(2) and the other give nothing new.

(v) Let $\Phi_1=A_3$ and $\Pi_1=\{-r_0, r_1, r_3\}$. Then $\Pi_2=\{r_2, r_5, r_6,r_7\}$ and $|N_W(\Pi_1)/W_2|=2$ by \cite[p. 197]{83Der}. Applying Lemma \ref{l:-1}, we conclude that $N_W(\Pi_1)=W_2\times \langle \theta(-1)\rangle$, so it is suffices to consider conjugacy classes of $W_2$. Also
$Z(\overline R)=\overline T_2$.

Using \cite[Corollary 1]{16But.t}, we see that every exponent in question divides $q^4-1$ or $(q^3\pm 1)(q\mp1)$. Since $W_2=W(A_1)\times W(A_3)\simeq S_2\times S_4$, every conjugacy class in $W_2$ corresponds to a~pair consisting of a partition of $2$ and a partition of $4$. If $H_w/\overline Z$ with $w\in W_2$ has exponent $q^4-1$, then it is clear that $w$
corresponds to the pair $[1,1]$, $[4]$. This case is handled  in the first part of Table \ref{tab:7a3}, and we see that the
exponent is in fact equal to $(q^4-1)/2$. On the other hand, the second part shows that the exponents $(q^3\pm 1)(q\mp1)$ do occur.

\tiny

\begin{table}

\caption{$\Phi_1=A_3$, $\Pi_1=\{-r_0,r_{1}, r_{3}\}$, $t_1=t_3=t_4=1$, $y=(1,1,0,1)$}
\label{tab:7a3}

$\begin{array}{|rclrcl|}
\hline
M_w&=&\begin{pmatrix}
1&0&0&0&0&0&0\\
0&1&0&0&0&0&0\\
0&0&1&0&0&0&0\\
0&0&0&1&1&0&0\\
0&0&0&0&0&1&0\\
0&0&0&0&0&0&1\\
0&0&0&0&-1&-1&-1\\
\end{pmatrix},&
X&=&\begin{pmatrix}
q - 1&0&0&0\\
0&-1&q&0\\
0&0&-1&q\\
0&-q&-q&-q - 1\\
\end{pmatrix},\\
A&=&\begin{pmatrix}
       1&  -q + 1&   -q + 1&        0\\
       0&        1&        0&        0\\
       0&        0&        1&        0\\
       0&       -q& -q^2 - q&        1\\
\end{pmatrix},
&
B&=&\begin{pmatrix}
    1  &  1&    1&    0\\
   0&   -1&   -q& -q^2\\
   0&    0&   -1&   -q\\
   0&    0&    0&   -1\\
\end{pmatrix},\\
y&\equiv& A[1]+A[4]\mypmod2,& yB&\equiv& (1,0,0,0)\mypmod2,\\
AXB&=&\diag(q-1,1,1,(q^2+1)(q+1)),& \exp&=&[2(q-1), (q^2+1)(q+1)]=(q^4-1)/2\\
\hline
M_w&=&\begin{pmatrix}
1&0&0&0&0&0&0\\
0&-1&0&0&0&0&0\\
0&0&1&0&0&0&0\\
0&1&0&1&0&0&0\\
0&0&0&0&1&1&0\\
0&0&0&0&0&0&1\\
0&0&0&0&0&-1&-1\\
\end{pmatrix},&
X&=&\begin{pmatrix}
-q - 1&0&0&0\\
0&q - 1&q&0\\
0&0&-1&q\\
0&0&-q&-q - 1
\end{pmatrix},\\
A&=&\multicolumn{4}{l|}{\begin{pmatrix}
      -q^2+q-1 &      q^3-q^2-q+2&       -q^2-1&       -1\\
      0 &      1&       0&       0\\
      0 &      0&       1&       0\\
\frac{q^3-1}{2} &  \frac{-q^4+2q^2+3q+2}{2}      &\frac{(q^2+1)(q+1)}{2}&   \frac{q+1}{2}\\
\end{pmatrix},}\\
B&=&\multicolumn{4}{l|}{\begin{pmatrix}
    1&      0 &     0&      \frac{q^3-1}{2}\\
     q^2& q^2-q-1&      -q&   \frac{q^5+q^2}{2}\\
    -q^2+q & -q^2+2q    & q - 1&   -\frac{(q^3+1)(q^2-q)}{2}\\
 -q + 1 &  -q+2&      1& -\frac{(q^3+1)(q-1)}{2}\\
\end{pmatrix},}\\
y&\equiv &A[1] \pmod 2,&yB&\equiv& (0,0,0,1)\mypmod2,\\
AXB&=&\diag(2,1,1,(q^3-1)(q+1)/2),& \exp&=&(q^3-1)(q+1)\\
\hline
\end{array}$
\end{table}
\normalsize

(vi) Let $\Phi_1=A_3^2$ and $\Pi_1=\{-r_0,r_1,r_3, r_5,r_6,r_7\}$. Then $\Pi_2=\{r_2\}$,  $Z(\overline R)$ consists of the elements of the form $h_2(t_2)h_5(t^3)h_6(t^2)h_7(t)$ with $t^4=1$ and $N_W(\Pi_1)\simeq \mathbb{Z}^2$ by \cite[p. 201]{83Der}.
By Lemma \ref{l:-1}, we know that $W_2\times \langle -w_0(\Pi_1)\rangle\leq N_W(\Pi_1)$. Also $w'=-w_0(\Pi_1\cap \{r_4\})$ lies in $N_W(\Pi_1)$. Thus we may assume that $w\in W_2\times \langle w'\rangle$. Also $w$ must permute irreducible components of $\Phi_1$, so $w=w'$ or
$w=w_{r_2}w'$.

Denoting $h_5(t^3)h_6(t^2)h_7(t)$ by $h(t)$, we verify that $h_2(t_2)h(t)w_{r_2}=h_2(t_2^{-1})h(t)$ and $h_2(t_2)h(t)w'=h_2(t_2t^2)h(t^{-1})$. So by \eqref{e:h}, it follows that $t_2^{\tau q-1}t_7^{2q}t_0=1$ and $t_7^{q-1}t_0=1$, where $\tau\in\{+,-\}$. If $q\equiv 1\pmod 4$, then the second equation yields $t_0=1$, and hence $H_w\simeq \mathbb{Z}_{2(q-\tau)}\times\mathbb{Z}_2$.
If $q\equiv -1\pmod 4$, then $t_0=t_7^2$ and the first equation yields $t_2^{\tau q-1}=1$, whence $H_w\simeq \mathbb{Z}_{q-\tau}\times \mathbb{Z}_4$.
In either case, the exponent of $H_w$ divides $2(q-\tau)$, which, in turn, divides $(q^4-1)/2$.

(vii) Let $\Phi_1=D_4$ and $\Pi_1=\{-r_0, r_1, r_3, r_9\}$. Then $\Phi_2=A_1^3$, $\Pi_2=\{r_2, r_6, r_{12}\}$, where $r_{12}=r_5+r_6+r_7$. According to \cite[p. 199]{83Der},
we have $N_W(\Pi_1)=W_2\rtimes S_3\simeq W(B_3)$, and using Lemma \ref{l:-1}, we see that $\theta(-1)$ corresponds  $-1\in W(B_3)$. Also $Z(\overline R)$ consists of the elements of the form $h_2(t_2)h_5(t_7)h_6(t_6)h_7(t_7)=h_2(t_2)h_6(t_6)h_{12}(t_7)$.

If $w\in W_2$, then it is clear that the exponent of $H_w$
divides $2[q-1,q+1]=q^2-1$, so $\eta(\Phi_1,w)$ divides $n_p(7)(q^2-1)$ from Item (5). There are six conjugacy classes outside  $W_2$, namely, those corresponding to  the signed cycle-types $[2,\overline 1]$, $[\overline 2,\overline 1]$, $[3]$,
$[2,1]$, $[\overline 2, 1]$, and $[\overline 3]$ (in notation introduced before Lemma \ref{l:dn}). It is sufficient to considered the first three types and the corresponding results are given in Table \ref{tab:7d4}. 

\tiny

\begin{table}

\caption{$\Phi_1=D_4$, $\Pi_1=\{-r_0,r_{1}, r_{3}, r_9\}$, $t_1=t_3=t_4=1$, $t_7=t_5$, $y=(1,1,0)$}
\label{tab:7d4}

$\begin{array}{|rclrcl|}
\hline
M_w&=&\begin{pmatrix}
 1&  0&  0&  0&  0&  0&  0\\
 0&  0&  0&  0&  0&  1&  0\\
 0&  1&  1&  2&  2&  1&  0\\
 0&  0&  0&  0&  0&  0&  1\\
 0& -1&  0& -1& -1& -1& -1\\
 0&  1&  0&  0&  0&  0&  0\\
 0&  0&  0&  1&  0&  0&  0\\
\end{pmatrix},&
X&=&\begin{pmatrix}
-1&0&q\\
-q&-q - 1&-q\\
q&0&-1\\
\end{pmatrix},\\
A&=&\begin{pmatrix}
1&0&0\\
1&1&0\\
q&0&1\\
\end{pmatrix},
&
B&=&\begin{pmatrix}
-1&0&q\\
1&-1&-q\\
0&0&1\\
\end{pmatrix},\\
y&\equiv& A[2]+A[3] \pmod 2,&yB&\equiv& (0,1,0)\mypmod2\\
AXB&=&\diag(1,q+1,q^2-1),& \exp&=&[2(q+1),q^2-1]=q^2-1\\
\hline
M_w&=&\begin{pmatrix}
1&0&0&0&0&0&0\\
0&0&0&0&0&1&0\\
0&1&1&2&2&1&0\\
0&0&0&0&0&0&1\\
0&0&0&-1&-1&-1&-1\\
0&-1&0&0&0&0&0\\
0&1&0&1&0&0&0\\
\end{pmatrix}, &
X&=&\begin{pmatrix}
-1&0&q\\
q&-q - 1&-q\\
-q&0&-1\\
\end{pmatrix},\\
A&=&\begin{pmatrix}
1 &  0 &  0\\
\frac{q^4-q^3+q^2+q-2}{4}&\frac{q^4+2q^2+1}{4}& \frac{q^4-q^3+q^2-3q-2}{4}\\
q &  q - 1 &  q - 2\\
\end{pmatrix},
&
B&=&\begin{pmatrix}
\frac{q^3 - q + 2}{2}&\frac{q^3-q}{2}&\frac{q^3 + q}{2}\\
 1&                   1&                   1\\
 \frac{q^2 - 1}{2}& \frac{q^2 - 1}{2}&\frac{q^2 + 1}{2}\\
\end{pmatrix},\\
y&\equiv& A[2]+A[3]\pmod 2,&yB&\equiv& (0,1,0)\mypmod2\\
AXB&=&\diag(-1,-(q^2+1)(q+1)/2,2)& \exp&=&(q^2+1)(q+1)\\
\hline
M_w&=&\begin{pmatrix}
  1&  0 & 0&  0& 0& 0& 0 \\
  0& 0& 0& 0& 0& 1& 0 \\
 -2&-2&-3&-4&-3&-2&-1 \\
  1& 2& 2& 3& 2& 1& 1 \\
  0&-1& 0&-1& 0& 0& 0 \\
  0& 0& 0& 0&-1&-1&-1 \\
  0& 0& 0& 1& 1& 1& 1 \\
\end{pmatrix},&
X&=&\begin{pmatrix}
-1&0&q\\
-q&q - 1&q\\
0&-q&-q - 1\\
\end{pmatrix},
\\
A&=&\begin{pmatrix}
1 & 0& 0 \\
q+1&-1&-1 \\
-q^4 - q^3 - q^2 + q + 1&q^3 + q - 1&q^3 + q - 2 \\
\end{pmatrix},&
B&=&\begin{pmatrix}
     -1 &    -q&     -q\\
     -1&q^2 + 2&q^2 + 1\\
      0 &    -1&     -1\\
\end{pmatrix},\\
y&\equiv& A[3]\mypmod2,& yB&\equiv& (0,0,1)\mypmod2,\\
AXB&=&\diag(1,1,q^3-1),& \exp&=&q^3-1\\
\hline

\end{array}$
\end{table}
\normalsize

(viii) Let $\Phi_1=(A_5)'$ and $\Pi_1=\{-r_0, r_1, r_2, r_3, r_4\}$. Then $\Pi_2=\{r_6, r_7\}$, $N_W(\Pi_1)=W_2\times \langle -w_0(\Pi_1)\rangle$  and
$Z(\overline R)=\overline T_2\times \overline Z$. It follows that $H_w/\overline Z$ is isomorphic to a maximal torus of $SL_3(q)$ or $SU_3(q)$, so
its exponent is equal to $q^3\pm 1$ or divides $q^2-1$.

(ix) Let $\Phi_1=(A_5)''$ and $\Pi_1=\{-r_0,r_1,r_3,r_4,r_5\}$. Then $\Pi_2=\{r_7\}$ and $N_W(\Pi_1)=W_2\times \langle -w_0(\Pi_1)\rangle$.  It suffices to consider the cases $w=1$ and $w=w_7$, and Table \ref{tab:7a5} shows that
all arising numbers divide $n_p(5)(q^2-1)$.

\tiny
\begin{table}

\caption{ $\Phi_1=(A_5)''$, $\Pi_1=\{-r_0,r_{1}, r_{3}, r_4, r_5\}$, $t_1=t_3=t_4=1$, $t_2=t_5^{-1}$, $t_6=t_5^2$, $y=(1,1)$}
\label{tab:7a5}

$\begin{array}{|rclrcl|}
\hline
M_w&=&\multicolumn{4}{l|}{I,\quad
X\ =\ \begin{pmatrix}
q-1&0\\
0&q-1\\
\end{pmatrix},\quad
A\ =\ \begin{pmatrix}
1& 0\\
1& 1\\
\end{pmatrix},\quad
B\ =\ \begin{pmatrix}
    1& 0\\
    -1&1\\
\end{pmatrix},}\\
y&\equiv&\multicolumn{4}{l|}{ A[2]\mypmod2,\quad yB\equiv (0,1)\mypmod2, \quad
AXB=X,\quad \exp=q-1}\\
\hline
M_w&=&\begin{pmatrix}
1&0&0&0&0&0&0\\
0&1&0&0&0&0&0\\
0&0&1&0&0&0&0\\
0&0&0&1&0&0&0\\
0&0&0&0&1&0&0\\
0&0&0&0&0&1&1\\
0&0&0&0&0&0&-1\\
\end{pmatrix},&
X&=&\begin{pmatrix}
q - 1&2q\\
0&-q - 1\\
\end{pmatrix},\\
A&=&\begin{pmatrix}
1& 1\\
0& 1\\
\end{pmatrix},
&
B&=&\begin{pmatrix}
1& 1\\
0&-1\\
\end{pmatrix},\\
y&\equiv& A[1]\pmod 2,&yB&\equiv& (1,0)\mypmod2,\\
AXB&=&\diag(q-1,q+1),& \exp&=&[q-1,q+1]=(q^2-1)/2\\
\hline
\end{array}$
\end{table}

\begin{table}

\caption{ $\Phi_1=D_5$, $\Pi_1=\{-r_0,r_{1}, r_{3}, r_4, r_9\}$, $t_1=t_3=t_4=1$, $t_5=t_2^{-1}$, $t_7=t_2^{-1}$, $y=(1,0)$}
\label{tab:7d5}

$\begin{array}{|rclrcl|}
\hline
M_w&=&\multicolumn{4}{l|}{I,\quad
X\ =\ \begin{pmatrix}
q-1&0\\
0&q-1\\
\end{pmatrix},\quad
A\ =\ I,\quad
B\ =\ I,}\\
y&\equiv&\multicolumn{4}{l|}{ A[1]\mypmod2,\quad yB\equiv (1,0)\mypmod2, \quad
AXB=X,\quad \exp=q-1}\\
\hline
M_w&=&\begin{pmatrix}
1&0&0&0&0&0&0\\
0&1&0&0&0&0&0\\
0&0&1&0&0&0&0\\
0&0&0&1&0&0&0\\
0&0&0&0&1&1&0\\
0&0&0&0&0&-1&0\\
0&0&0&0&0&1&1\\
\end{pmatrix},&
X&=&\begin{pmatrix}
q - 1&-2q\\
0&-q - 1\\
\end{pmatrix},\\
A&=&\begin{pmatrix}
1& 0\\
\frac{-q-1}{2}& 1\\
\end{pmatrix},
&
B&=&\begin{pmatrix}
2& -q\\
1&\frac{-q+1}{2}\\
\end{pmatrix},\\
y&\equiv& A[1]\pmod 2,&yB&\equiv& (0,1)\mypmod2,\\
AXB&=&\diag(-2,(q^2-1)/2),& \exp&=&q^2-1\\
\hline
\end{array}$
\end{table}
\normalsize

(x) Similarly, if $\Phi_1=D_5$ and $\Pi_1=\{-r_0, r_1, r_3, r_4, r_9\}$, then $\Pi_2=\{r_6\}$ and  $N_W(\Pi_1)=\langle w_6\rangle\times \langle -w_0(\Pi_1)\rangle$, and we may assume that $w=1$ or $w=w_6$. By Table \ref{tab:7d5},
we see that $n_p(7)(q^2-1)\in\omega(G)$.

(xi) If $\Phi_1=A_7$, then $Z(\overline R)$ is a cyclic group of order 4, so $\eta(A_7,w)$ divides $n_p(7)(q^2-1)$.

(xii) Let $\Phi_1=D_6$ and $\Pi_1=\{-r_0,r_1,r_2,r_3,r_4,r_5\}$. Then $\Pi_2=\{r_7\}$, $N_W(\Pi_1)=W_2$ and $Z(\overline R)=\overline T_2\times \overline Z$. We conclude that $\eta(D_6,w)$ is equal to $n_p(7)(q\pm 1)$, which  divide $n_p(11)(q\pm 1)$.

(xiii) Let $\Phi_1=E_6$ and $\Pi_1=\{-r_0,r_1,r_3,r_4,r_5, r_8\}$. Then $\Phi_2$ is empty and $N_W(\Pi_1)=\langle -w_0(\Pi_1)\rangle$, so we may assume that $w=1$. Since $Z(\overline R)$ consists of elements of the form $h_2(t^{-1})h_5(t)h_6(t^2)h_7(t)$, we see that $H_w$ is a cyclic group of order $2(q\mp 1)$. Together with the fact that $mh(E_6)=11$, this gives Item (6).

(xiv) If $\Phi_1=E_7$, then $\eta(\Phi_1,w)=n_p(17)$, so this provides Item (7) and completes the proof of the first statement of the theorem.

\subsection{The set $\omega(G\setminus L)$} To prove that $\nu_\delta(q)\subseteq\omega(G\setminus L)$, we need the following lemma.

\begin{lemma}\label{l:diag}
Suppose that $\overline Z\leqslant Z(\overline R)^0$. Then $\eta(\Phi_1,w)\in\omega(G\setminus L)$ for every $w\in N_W(W_1)$.
\end{lemma}

\begin{proof}
Denote $Z(\widetilde R)_{\sigma\circ w}$ by $H$ and $Z(\overline  R)_{\sigma\circ w}/\overline Z$ by $H_0$. It is clear that $H_0=H\cap L$, and so $|H:H_0|\leqslant 2$. To prove the lemma, it suffices to show that $\exp(H)\in\omega(H\setminus H_0)$. Since $z\in Z(\overline R)^0$, by the Lang--Steinberg theorem \cite[Theorem 10.1]{68Ste}, there is $g\in Z(\overline R)^0$ such that $g^{\sigma\circ w}g^{-1}=z$.
It follows that $g\in H\setminus H_0$ and thus $|H:H_0|=2$. Now it is easy to see that there is an element $h\in H\setminus H_0$ of order $\exp(H)$.
\end{proof}

Now we observe that every number of $\nu_\delta(q)$ arises in a reductive group with connected center: we work in tori for the numbers in Item (1) and we can take $\Phi_1$ to be equal to $A_1$, $A_3$, $D_4$, $D_5$, $E_6$ for the numbers in Items (2)--(6) respectively, with Table \ref{tab:CondE7} showing that the center is connected in these cases. Thus $\nu_\delta(q)\subseteq \omega(G\setminus L)$.

If $a\in\omega(G\setminus L)$, then $a$ is even and so divides some even number of $\nu(q)$.
But $\nu_\delta(q)$ is exactly the subset of even elements of $\nu(q)$, and therefore $\mu(\omega(G\setminus L))\subseteq \nu_\delta(q)$. This completes the proof of Theorem \ref{t:E7}.

\section{Proof of Theorem \ref{t:main}}

We begin with a proof of Proposition \ref{p:f}. The first claim of the proposition is a corollary of \cite[Proposition 13]{06Zav.t}, so we only need to prove the second one.

Suppose that $q=p^m$ with odd prime $p$. Recall from Section \ref{s:diag} that $L=\overline L_\sigma/\overline Z$ and $\Inndiag L=(\overline L/\overline Z)_\sigma$. Denote by $K$ the preimage of $(\overline L/\overline Z)_\sigma$ in $\overline L$. Also recall from Introduction that $\varphi$ is the field automorphism of $L$ mapping $x_r(s)$ to $x_r(s^p)$. It is clear that $\varphi$ is induced from the corresponding endomorphism of $\overline L$, and in what follows $\varphi$ denotes this endomorphism, and so $\sigma=\varphi^m$.

Now let $\psi\in\langle \varphi\rangle$ and $g\overline Z\in\Inndiag L\setminus L$. Then $g\in K\setminus\overline L_\sigma$ and so $g^{\sigma}g^{-1}=z$, where $z$ is the generator of $\overline Z$. We may assume that $\psi^k=\sigma$, where $k$ is the order of $\psi$ regarded as an automorphism of $L$.
By the Lang--Steinberg theorem \cite[Theorem 10.1]{68Ste}, there is $y\in \overline L$ such that $g=y^{-\psi}y$, and we define $h=yy^{-\sigma}$. Then \begin{equation}\label{e:order}
(\psi g\overline Z)^k=g^{\psi^{k-1}}\dots g\overline Z= y^{-\psi^k}y^{\psi^{k-1}}y^{-\psi^{k-1}}y^{\psi^{k-2}}\dots y^{-\psi}y\overline Z=y^{-\sigma}y\overline Z=y^{-1}h\overline Zy.\end{equation}
It follows that $|\psi g\overline Z|=|\psi|\cdot |h\overline Z|$. Also we have \begin{equation} \label{e:fix} g^\sigma =gz\Leftrightarrow y^{-\psi\sigma}y^{\sigma}=y^{-\psi}yz\Leftrightarrow y^\psi y^{-\psi\sigma}=yy^{-\sigma}z \Leftrightarrow h^\psi=hz.\end{equation}
Therefore, $h\overline Z\in (\overline L/\overline Z)_\psi$ and $h\not\in\overline L_\psi$, which yields $|h\overline Z|\in\omega(\Inndiag L_0\setminus L_0)$. Thus
$\omega(\psi (\Inndiag L\setminus L))\subseteq |\psi|\cdot \omega(\Inndiag L_0\setminus L_0)$.

And vice versa, if we have $h\in\overline L$ with
$h^\psi h^{-1}=z$, choose $y\in \overline L$ such that $h=yy^{-\sigma}$ and define $g=y^{-\psi}y$, then \eqref{e:order} and \eqref{e:fix} guarantee that $g\in K\setminus \overline L_\sigma$ and $|\psi g\overline Z|=|\psi|\cdot |h\overline Z|$. This gives the reverse containment, and completes the proof of Proposition \ref{p:f}.

By Proposition \ref{p:f}, to prove Theorem \ref{t:main}, it suffices to
consider the cases $\alpha\in L$ and $\alpha\in\delta L$. These cases are covered by  \cite[Theorem 2]{16But.t} and the second assertion of Theorem~\ref{t:E7} respectively.

\end{document}